\documentclass[a4paper, 11pt]{amsart}

\usepackage{geometry}

\newcommand{\dsum}{\displaystyle\sum}

\usepackage{tikz}
\usetikzlibrary{shapes,arrows}
\usepackage{rotating}
\usepackage{multirow}
\usepackage{hyperref}
\usepackage{booktabs}

\usepackage{amssymb}
\usepackage{amsmath,amsthm}
\usepackage{amsfonts}
\let\origmaketitle\maketitle
\def\maketitle{
  \begingroup
  \def\uppercasenonmath##1{} 
  \let\MakeUppercase\relax 
  \origmaketitle
  \endgroup
}
\begin{document}

		\title[A selection system based on player performance]{\large A multicriteria selection system based on player performance. Case study: The Spanish ACB Basketball League}

		\author[V\'ictor Blanco, Rom\'an Salmer\'on \MakeLowercase{and} Samuel G\'omez-Haro]{{\large V\'ictor Blanco$^\dagger$, Rom\'an Salmer\'on$^\dagger$ \MakeLowercase{and} Samuel G\'omez-Haro$^\ddagger$}\medskip\\
$^\dagger$Dpt. Quant. Methods for Economics \& Business, Universidad de Granada, SPAIN.\\
$^\ddagger$Dpt. of Business \& Management, Universidad de Granada,
}

\address{V. Blanco: Dpt. Quant. Methods for Economics \& Business, Universidad de Granada}
\email{vblanco@ugr.es}
\address{R. Salmer\'on: Dpt. Quant. Methods for Economics \& Business, Universidad de Granada}
\email{romansg@ugr.es}
\address{S. G\'omez-Haro: Dpt. of Business \& Management, Universidad de Granada}
\email{vblanco@ugr.es}

\begin{abstract}
In this paper, we describe an approach to rank sport players based on their efficiency. Although is extremely useful to analyze the performance of team games there is no unanimity on the use of a single index to perform such a ranking. We propose a method to summarize the huge amount of information collected at different aspects of a sport team  which is almost daily publicly available. The tool will allow agents involved in a player's negotiation to show the strengths (and weaknesses) of the player with respect to other players. The approach is based on applying a multicriteria outranking methodology using as alternatives the potential players and criteria different efficiency indices. A novel automatic parameter tuning approach is detailed that will allow coaches and sports managers to design templates and sports strategies that improve the efficiency of their teams. We report the results performed over the available information on the ACB Basketball League, and we show how it can be easily implemented and interpreted in practice by decision-makers non familiar with the mathematical side of the methodology.

\end{abstract}

\maketitle

\section{Introduction}
\label{intro}

Operations Research is one of the most active fields to help decision makers with objective tools for choosing the ``best'' of their alternatives based on the evaluation of the adequate criteria. A wide family of tools are available in the literature to make the best decision, each of them depending of the type of provided information, characteristic of the decision, number of agents involved on it or the number of criteria and alternatives, etc. Needless to say that in the current digital age in which the amount of available data is exponentially increasing, it is specially important to provide quantitative tools to allow practitioners an efficient management of such an information. In particular, in Sport Science, each team or individual player updates, almost daily, its available information. Furthermore, making the best decision instead of a nearly good may induce a significative economic or social gain to the teams, player agents, or any of the several agents involved in the Sport Industry. Hence, using adequate tools to summarize and manage the information  is one of the most important tasks in Sport Management. Multicriteria Analysis allows one to consider at the same time several criteria to chose among the available alternatives. If the number of alternatives is of reasonable size  several tools are available from the so-called discrete multicriteria/multiattribute analysis, that allow one to rank, at least partially, the alternatives (see \cite{electre} and \cite{ahp}, among many others). Observe that, unless a preference weight is provided for each of the criteria that allows to unify the results of all of them into a single utility value, for a given alternative, the best result for a criterion may not coincide with the best choice for other criteria.

In particular, the goal of this paper is to adapt one the available ranking tools to measure the performance of players, based on the available criteria, as the annotated points, playing time or number of tries, and many others. Agents of basketball players highlighting the value of its players when negotiating with club presidents, would be helped with quantitative measures about the player on its current team (see \cite{j2016}).  Also, coaches and managers of sports teams design the lineup of players, and at the end of each season they evaluate the continuation of the players  as well as the incorporation of new talents to complement and improve the team. Even for short-terms, the coaches select the most suitable players for a game or part of it. In order to make the decisions, among other factors, they are provided with a lot of  information and the correct decision is expected by analyzing the available data. A deeper understanding of resources implies an adequate decision and then an improvement of the performance of the team (see \cite{Alamar} and Drust\cite{Drust}), since it allows one to detect the strengths and weaknesses of the players and the whole team.

One of the main differences between basketball teams and other sports comes from the availability of information. One can find data about each individual player or team for a game,  from the number of scored points to the completed possessions. Hence, the detection of factors implying success in a basketball game is particularly hard since it depends on the accomplishment of numerous game actions (see \cite{Puenteetal2015}). Hence, in this framework, it is crucial to analyze the performance of teams and players by summarizing the available information. Thus, several desirable skills for players selection should be adequately selected to be summarized. However, although different indicators have been proposed in the literature, there is not an unified measure for the performance of a basketball player, the most popular being the TENDEX index \cite{Heeren}, the ratio of performance-performance index rating used in competitions of European basketball (roughly, the difference between positive and negative action games), the Player Efficiency Rating the (PER) of Hollinger, or different valuation plus/minus rates \cite{Oliver}. Note that an index to measure the performance of a player must take into account the characteristics of the player, its  experience in the game, its profile (offensive or defensive), etc (see \cite{Gutierrez}), and then, an unified indicator is not possible for the whole set of players. Some discussions about the suitability of each of the proposed measures has been partly addressed in Berri and Bradbury\cite{Berri}, Lewin and Rosenbaum\cite{Lewin} and Winston\cite{Winston}, among others.

Although the selection of the adequate or more representative indices is highly interesting in sport science, as stated in \cite{Mertzetal2016}, computing an global efficiency index is difficult due to the variety of individual statistics. In \cite{Zivetal2010}, the authors recommen to limit the use of individual statistics for predicting the final ranking of their teams. We propose to turn this disadvantage into a strength by considering a different perspective to provide a preference ranking of the players based on using all the available information. More concretely, we use tools borrowed from Multicriteria Analysis and adapt them to help managers and/or coaches in the players' selection task or for agents of player that wishing to highlight their clients' strengths over other players. As far as we know, these tools have not been applied in this way to analyze sport players.

In particular, we adapt the PROMETHEE method, introduced by Brans~\cite{Brans0}. The PROMETHEE method have been successfully applied in several fields, as energy sustainability \cite{m2016} or Education \cite{z2016}. It allows one to aggregate several criteria, whose results are known for each given player and season, into one or two indices, the so-called \textit{net flow preference indices}. With this method we provide a support tool for practitioners which are not familiar which the quantitative techniques the approach, but that can be intuitively used.  Hence, describe a method for selecting players which complement the usual techniques in this framework. One of the main highlights of our approach is that we propose a simple automatic strategy to determine the parameters needed for the application of the PROMETHEE methodology: the thresholds and the weights. This will allow agents and coaches which are not familiar with the quantitative tools to easily use and interpret the results derived from the application of the approach.

We apply the methodology to the case study of the 191 players who participated in the Spanish ACB Basketball League during the 2014-2015 season. We construct a ranking system by classifying the players into five different profiles (positional classification) and analyze two different scenarios for each of them. Different graphical representations of the results are shown which allow for a better understanding of the insights of our approach.

The paper is organized as follows. In Section \ref{sec:promethee} we recall the elements and philosophy of the PROMETHEE methodology. Section \ref{sec:casestudy} is devoted to apply the multicriteria methodology to our real-world case study. Finally, in Section \ref{sec:conclusions} we draw some conclusions and further research on the topic.

\section{The outranking methodology}
\label{sec:promethee}

MultiCriteria Decision Aid (MCDA) allows one ranking a given finite set of alternatives based on multiple (probably) conflicting criteria. Several MCDA approaches have been proposed in the literature (see, for instance, \cite{mcda}) and they have been applied to real-world problems in many different areas, as management or engineering, amongst others. The main problem and motivation, in the development of methods where different criteria are involved in the decision process, is that the best in one criterion may imply a worsering in other,. Hence, one is not able to determine which is the best or worst alternative. In general, only  a partial order is possible with such an information. One of the choices in this framework is to provide a set of weights of relative importance for each criterion, and then, if the measures are of quantitative nature, one may compute an overall index to summarize the given criteria. However, this is in general not possible, and when it is, it is not easy to implement. One of the methods designed to rank a set of alternatives based on different criteria is the PROMETHEE method (Preference Ranking Organization Method for Enrichment Evaluations), that was introduced by Brans~\cite{Brans0}. There, two first versions of the method are provided, namely, PROMETHEE I and PROMETHEE II. Several extensions have appeared since then, from PROMETHEE III, for ranking based on interval, to PROMETHEE CLUSTER for nominal classification, in a series of paper mostly coauthored by Brans and Mareschal\cite[amongst other]{BM1,BM2,BM3}. Also, there has appeared some recent developments related to the PROMETHEE method as in the one by Soylu\cite{soylu}, in which a Tchebycheff agreggation technique is applied, the survey by Saaty and Ergu\cite{saaty-ergu} where several multicriteria methods are compared in order to make a \textit{right} decision or the combination of the SMAA and PROMETHEE to determine adequate weights and thresholds \cite{corrente}, to mention a few. The importance of the method in the Operations Research community is also reflected in that there is an available software, Visual PROMETHEE\cite{VP}, that allows one to perform all the computations. Actually, the flexibility and applicability in many contexts of PROMETHEE methodologies have made it more and more popular in the recent years (see Brans and Mareschal\cite{Brans} and the references therein for a complete overview of the method).

To be self-contained we detail the general performance of PROMETHEE methods. Let us consider a finite set of alternatives $\mathcal{A}=\{a_1, \ldots, a_n\}$ and a set of criteria to be, without loss of generality, maximized to make a decision over the set of alternatives, $\mathcal{C}=\{c_1, \ldots, c_m\}$.  Every alternative $a_i$ is assumed to be evaluated under each of the criteria $c_k$, with value $r_{ik}$, for $i=1, \ldots, n$, $k=1, \ldots, m$.  The PROMETHEE paradigm is based on the principle that the decision maker may manifest some subjective information over the criteria that allows for evaluating the results under its own perspective. Each pair of alternative is then compared for each of the criteria, computing the absolute difference on their results:
$$
d_{ij}^k = r_{ik}-r_{jk}.
$$
Those differences are now evaluated using one of the six classical preference functions which are detailed in Table \ref{criteria}. A preference function for the $k$-th criterion is a mapping $H_k: \mathbb{R} \rightarrow [0,1]$ that measures the intensity of agreement with the statement $a_i$ is better than $a_j$ according to criterion $c_k$. A $H_k$-value of $0$ will imply that the alternatives are indifferent while if the $H_k$-value is $1$, it will imply that one of the alternatives (the one with larger result) is highly preferred to the other. Intermediate values in $(0,1)$ determine the degree of preference between the two alternatives.

\begin{table}
\centering
\caption{Preference Criteria for PROMETHEE method}   \label{criteria}
\begin{tabular}{ccc}\toprule
Type I:  Usual &
$H (d) = \left\{\begin{array}{rl}
0 & \mbox{ if $d=0$}\\
1 & \mbox{ otherwise}\end{array}\right.
$ &
\begin{minipage}[c]{3cm}
\begin{center}
\begin{tikzpicture}[scale=0.65]

\draw[->, thin] (0,0)--(2,0);
\draw[<-, thin] (0,1.5)--(0,0);
\draw[thin] (-0.05,1)--(0.05,1);
\draw[gray, thick] (0.05,1)--(2,1);
\draw[gray] (0,1) circle (2pt);
\fill[gray] (0,0) circle (2pt);

\node[left] at (0,0) {\tiny $0$};
\node[left] at (0,1) {\tiny $1$};
\node[right] at (2,0) {\tiny $d$};
\end{tikzpicture}
\end{center}
\end{minipage}\\
Type II: U-shape &
$H (d) = \left\{\begin{array}{rl}
0 & \mbox{ if $d\leq q$}\\
1 & \mbox{ otherwise}\end{array}\right.
$
 &\begin{minipage}[c]{3cm}
\begin{center}
\begin{tikzpicture}[scale=0.65]

\draw[->, thin] (0,0)--(2,0);
\draw[<-, thin] (0,1.5)--(0,0);
\draw[thin] (-0.05,1)--(0.05,1);
\draw[gray, thick] (1.05,1)--(2,1);
\draw[gray, thick] (0,0)--(1,0);
\fill[gray] (1,0) circle (2pt);
\draw[gray] (1,1) circle (2pt);

\node[left] at (0,0) {\tiny $0$};
\node[left] at (0,1) {\tiny $1$};
\node[right] at (2,0) {\tiny $d$};
\node[below] at (1,0) {\tiny $q$};
\end{tikzpicture}
\end{center}
\end{minipage}
\\
Type III: V-shape &
$H (d) = \left\{\begin{array}{rl}
\dfrac{d}{p} & \mbox{ if $d\leq p$}\\
1 & \mbox{ otherwise}\end{array}\right.
$  &
 \begin{minipage}[c]{3cm}
\begin{center}
\begin{tikzpicture}[scale=0.65]

\draw[->, thin] (0,0)--(2,0);
\draw[<-, thin] (0,1.5)--(0,0);
\draw[thin] (-0.05,1)--(0.05,1);
\draw[thin] (1,-0.05)--(1,0.05);
\draw[gray, thick] (1,1)--(2,1);
\draw[gray, thick] (0,0)--(1,1);

\node[left] at (0,0) {\tiny $0$};
\node[left] at (0,1) {\tiny $1$};
\node[right] at (2,0) {\tiny $d$};
\node[below] at (1,0) {\tiny $p$};
\end{tikzpicture}
\end{center}
\end{minipage}\\
Type IV: Level & $H (d) = \left\{\begin{array}{rl}
0 & \mbox{ if $d\leq q$}\\
0.5 & \mbox{ if $q < d\leq p$}\\
1 & \mbox{ otherwise}\end{array}\right.
$  &
\begin{minipage}[c]{3cm}
\begin{center}
\begin{tikzpicture}[scale=0.65]

\draw[->, thin] (0,0)--(2,0);
\draw[<-, thin] (0,1.5)--(0,0);
\draw[thin] (0.75,-0.05)--(0.75,0.05);
\draw[thin] (1.5,-0.05)--(1.5,0.05);
\draw[thin] (-0.05,1)--(0.05,1);
\draw[thin] (-0.05,0.5)--(0.05,0.5);
\draw[gray, thick] (0,0)--(0.75,0);
\draw[gray, thick] (0.8,0.5)--(1.5,0.5);
\draw[gray, thick] (1.55,1)--(2,1);
\fill[gray] (0.75,0) circle (2pt);
\fill[gray] (1.5,0.5) circle (2pt);

\draw[gray] (0.75,.5) circle (2pt);
\draw[gray] (1.5,1) circle (2pt);
\node[left] at (0,0) {\tiny $0$};
\node[left] at (0,0.5) {\tiny $0.5$};
\node[left] at (0,1) {\tiny $1$};
\node[right] at (2,0) {\tiny $d$};
\node[below] at (0.75,0) {\tiny $q$};
\node[below] at (1.5,0) {\tiny $p$};
\end{tikzpicture}
\end{center}
\end{minipage}\\
Type V: V-shape with Indifference &
$H (d) = \left\{\begin{array}{rl}
0 & \mbox{ if $|d|\leq q$}\\
\frac{d-q}{p-q} & \mbox{ if $q < d\leq p$}\\
1 & \mbox{otherwise}\end{array}\right.
$  &
\begin{minipage}[c]{3cm}
\begin{center}
\begin{tikzpicture}[scale=0.65]

\draw[->, thin] (0,0)--(2,0);
\draw[<-, thin] (0,1.5)--(0,0);
\draw[thin] (-0.05,1)--(0.05,1);
\draw[thin] (0.75,-0.05)--(0.75,0.05);
\draw[thin] (1.5,-0.05)--(1.5,0.05);
\draw[gray, thick] (0,0)--(0.75,0);
\draw[gray, thick] (0.75,0)--(1.5,1);
\draw[gray, thick] (1.5,1)--(2,1);

\node[left] at (0,0) {\tiny $0$};

\node[left] at (0,1) {\tiny $1$};
\node[right] at (2,0) {\tiny $d$};
\node[below] at (0.75,0) {\tiny $q$};
\node[below] at (1.5,0) {\tiny $p$};
\end{tikzpicture}
\end{center}
\end{minipage}
\\
Type VI: Gaussian &
$
H (d) = 1 - {\rm e}^{-\frac{d^2}{2\sigma^2}}
$ &
 \begin{minipage}[c]{3cm}
\begin{center}
\begin{tikzpicture}[scale=0.65]

\draw[->, thin] (0,0)--(2,0);
\draw[<-, thin] (0,1.5)--(0,0);
\draw[thin] (-0.05,1)--(0.05,1);
\draw[domain=-0:2,smooth,variable=\x,gray,thick] plot ({\x},{1 - exp(-\x*\x/2)});

\node[left] at (0,1) {\tiny $1$};
\node[right] at (2,0) {\tiny $d$};
   \end{tikzpicture}

\end{center}
\end{minipage}\\ \bottomrule
\end{tabular}
\end{table}

Observe that each of the preference criteria depends on some thresholds that allows the decision maker to incorporate its preferences. In particular, the $q$-parameter, usually called \textit{indifference threshold}, indicates an upper bound (in terms of the difference of results for a pair of alternative under a criterion) for indifference between alternatives. On the other hand, the $p$-parameter, so-called the \textit{preference threshold}, is a lower bound to model whether one of the two alternatives is clearly preferable to a second one under the given criterion. The most popular type of preference function for quantitative information is Type V (V-shape with indifference criterion). This function is parameterized with two values: the indifference ($q$) and the preference ($p$) threshold. The correct specification of the set of thresholds from which two alternatives are considered equal or different is the other important decision faced by the decision maker with this methodology. Different alternatives for fixing the threshold values are possible. One of them is based on the decision-maker experience, who set the values using its own subjective opinion. Other, is to provide a mechanism to automatically compute them from using the given sample. In what follows we describe a methodology to specify the thresholds by using some descriptive indices from the particular instance under the study.

Let us consider one criterion, $c_k$, which is assumed to be evaluated for all the alternatives. Let $D_k$ be the set of differences between all distinct alternatives for that criterion. We define the indifference and preference parameters as follows:
$$
q_k = Q_{\alpha}(D_k) \quad p_k = Q_{\beta}(D_k),
$$
where $Q_z$ denotes the $z\%$ quantile, for $z \in [0,100]$.

Note that this choice assures, among all the possible differences of the values of the criterion for all the alternatives (roughly, $n^2$), $(\beta-\alpha)\%$ of them are classified as not indifferent nor prioritary, while $\alpha\%$ and $(1-\beta)\%$ of the differences are considered as indifferent and preferential, respectively. This settings allows one to control that $(\beta-\alpha)\%$ of combinations of two players are in the \textit{positive slope} side of the Type V graph (see Table \ref{criteria}).

Once a preference criteria is specified for each criterion, $H_1, \ldots, H_m$, with their respective thresholds, for each of the criteria in the decision process, they are applied to the previously computed differences. Thus, the absolute differences are translated into preference values in $[0,1]$. Then, with those values, for each pair of alternatives a preference index is computed as follows:
$$
p_{ij} = \dsum_{k \in \{1, \ldots, m\}} \omega_k H_k(d_{ij}^k), \quad \forall i, j=1, \ldots, n,
$$
where $\omega_1, \ldots, \omega_m \geq 0$ represent the importance level of each criterion. It is commonly assumed that $\dsum_{k=1}^m \omega_k = 1$ to represent proportionality of importance of the criteria in the decision making process.

The preference flows are computed as:
$$
\phi^+ (a_i) = \dfrac{1}{n-1} \dsum_{j=1}^n p_{ij} \quad \mbox{and} \quad \phi^- (a_i) = \dfrac{1}{n-1} \dsum_{j=1}^n p_{ji},  \quad \forall i=1, \ldots, n.
$$
Observe that $\phi^+(a_i)$ is the average level of preference of the alternative $a_i$ with the given decision maker specifications, with respect to the rest of the alternatives. On the other hand, $\phi^-(a_i)$ is the average level of preference of the rest of the alternatives with respect to $a_i$. With those values, the PROMETHEE I method identifies if an alternative $a_i$ is preferable to $a_j$ if and only if $\phi^+(a_i) \geq \phi^+(a_j)$ and $\phi^-(a_i) \leq \phi^-(a_j)$ (with at least one of the inequalities being strict). This first method allows us to partially rank the alternatives. Note that if one of the two conditions does not hold, the alternative are not comparable with this approach.

Finally, the PROMETHEE II method uses the so-called net flow index:
$$
\phi(a_i) = \phi^+(a_i) - \phi^-(a_i),  \quad \forall i=1, \ldots, n.
$$
With such a flow, it is said that $a_i$ if preferable to $a_j$ if and only if $\phi(a_i) > \phi(a_j)$.

Finally, we summarize the elements needed to applied the proposed approach for ranking players based on the efficiency. A particular choice of the input information will be described in our case study.
\begin{itemize}
\item {\it Alternatives}: Players of a team during one or more seasons. In many sports is convenient to split players according to their positions, since the performance measure differ on their importance in the final ranking decision.
\item {\it Criteria}: For each of the alternatives (players), we must be able to measure two or more criteria which make impact on the efficiency of the player. An adequate choice of the criteria is crucial for a satisfying result, since they all put together conform the ranking system.  They are usually chosen by experts (as agents or coaches) which has some evidence about which measures of efficiency are crucial in determining the a global efficiency measure, but also they must be in concordance with the available information.
\item {\it Aggregation weights}: Other aspect where the decision-maker may incorporate its expertise to the methodology is on the establishment of the importance weights of each of the considered criteria. A very simple choice is to consider that all criteria are equally important. However, it is very usual to consider that some criteria are more relevant than others, with respect to the impact in the global efficiency. This is the case when the decision maker's goal is to rank a particular profile of alternatives, in which the same criteria are considered, but some of them are more important than other for the different play positions.
\item {\it Shape of the preference functions}: The decision maker should be able to choose one of the criteria in Table \ref{criteria}, to measure the degree of preference/indifference between players. The determination of the thresholds can be performed using the approach described in Section \ref{sec:promethee}, just by fixing the values of $\alpha$ and $\beta$. Other options are also possible, as the use of the SMAA method described in \cite{corrente}.
\end{itemize}

Using the above information, provided by the decision maker, one can apply the described methodology to obtain a partial (PROMETHEE I) or total (PROMETHEE II) ranking system. In order to illustrate the proposed approach, in the next section, we describe its application for ranking basketball players of different typologies.

\section{Case Study: Spanish ACB Basketball Players}
    \label{sec:casestudy}

In this section we apply the method described above to rank basketball players of different typologies based on several criteria. All the data used for the case study is available upon request for the interested readers.

The input data for this case study are the following:
\begin{itemize}
\item {\it Alternatives}: We consider as alternatives each of the 191 players who participated in the Spanish ACB Basketball League during the 2014-2015 season. For the sake of data consistency and availability, we do not consider players that did not played at least 10 games and an average of at least 10 minutes. These requirements are usually assumed in the literature for analyzing the quality of basketball players (see Berri et al Berri \textit{et. al.}\cite{Berri2007}, Berri and Krautmann\cite{Berri2006}, Cooper \textit{et. al.}\cite{Cooper2009}). The data were obtained from the official database of the Spanish ACB Basketball League (\url{http://www.acb.com}). In order to compare and adequately weigh the criteria, the players were split based on their play positions into five categories: point guards (42), shooting guards (38), forwards (33), power-forwards (36) and centers (32). 
\item {\it Criteria}: A set of 6 evaluation criteria have been considered in the case study for each player.  We use the following indices as criteria for the proposed methodology, which, as far as we know, have not been previously used:
\begin{description}
    \item[\bf PtsM: ]	Ratio of points scored by the player with respect to the number of minutes played: $$PtsM = \frac{Pts}{Min},$$
        where ${\rm Pts}$ and ${\rm Min}$ are, respectively, the points and minutes played by the player.
    \item[\bf DRM: ] Ratio of the player defensive rating  with respect to the number of minutes played: $${\rm DRM} = \frac{\rm DRB + STL + BLK - PF}{\rm Min},$$
        where ${\rm DRB}$, $STL$, $BLK$ and $PF$ are, respectively, the defensive rebounds, steals, blocks and personal fouls made by the player.
    \item[\bf ORM: ]	Ratio of the player offensive rating  with respect to the number of minutes played:
        \begin{eqnarray}
            {\rm ORM} &=& \frac{2 \cdot 2{\rm P} + 3 \cdot 3{\rm P + FT + ORB + AST + PFR}}{\rm Min} \nonumber \\
                & &  - \frac{\rm (FGA - FG) + (FTA - FT) + TOV + BLK}{\rm Min}, \nonumber
        \end{eqnarray}
        where $2P$, $3P$, $FT$, $ORB$, $AST$, $PFR$, $FGA$, $FG$, $FTA$, $TOV$ and $BLKR$ are, respectively, the 2-points field goals, 3-point field goals, free throws, offensive rebounds, assists, personal fouls received, field goal attempts, field goals (includes both 2-point field goals and 3-point field goals), free throw attempts, turnovers and blocks received by the player.
        \item[\bf EPts: ] Points efficiency measured as the ratio of points scored by the player with respect to the points that could have scored depending on the shots made and multiplied by 100: $$EPts = \frac{Pts}{2\cdot 2PA + 3 \cdot 3PA + FTA} \cdot 100,$$
        where $2PA$ and $3PA$ are, respectively, the 2-point field goal attempts and 3-point field goal attempts made by the player.
    \item[\bf ASTM: ]	Ratio of the number of assists and steals for every turnovers with respect to the number of minutes played: $${\rm ASTM} = \frac{\rm AST + STL}{{\rm TOV} \cdot {\rm Min}}.$$
    \item[\bf PCS \%: ] Ratio between the number of possessions completed successfully (PCS, when a player makes a field goal, receives a personal foul or assists) and the number of completed possessions (PC, when a player attempts a field goal, receives a personal foul, assists or lose a ball) multiplied by 100: $${\rm PCS} \% = \frac{\rm PCS}{\rm PC} \cdot 100 = \frac{\rm FG + PFR + AST}{\rm FGA + PFR + AST + TOV} \cdot 100.$$
\end{description}
Note that we split players by their game position since each profile has very particular characteristics, implying a significative difference in some of the measured criteria. In this way, we avoid bias between players due to the used criteria. For instance, point guard players usually get higher values in the ASTM index, while pivots obtain higher results in the DRM index, because its game position. This assertion was statistically checked and shown in Table \ref{anova}, where we compute the average of each criterion for each player's profile and we apply the ANOVA test to those values. The results indicates that means are statistically different for criteria DRM, ORM, EPts, ASTM and PCS \%.

\item {\it Aggregation weights}: For illustrative purposes, we consider the following two different scenarios in our case study:
\begin{itemize}
    \item Scenario 1: All the criteria are identically weighted.
    \item Scenario 2: Two criteria are overweighted for each of the players position clusters with the following strategy. From the difference between scored and received points by his team while it is on court (plus/minus, PM) weighted by minutes played divided by duration of the game, normally 40 minutes: $${\rm PMW} = {\rm PM} \cdot \frac{\rm Min}{40},$$ and considering that ${\rm PMW}$ is positively affect to the player value for all players profiles, we will consider that the above criteria that are positively related to $PMW$ must be higher weight. In Table \ref{table1} are shown the correlations between all criteria and ${\rm PMW}$ by player position  (the correlations significantly different from zero to the level of significance of 5\% are boldfaced). As can be observed, high values for ${\rm EPts}$ and ${\rm ASTM}$ (for exterior players) and high values for ${\rm DRM}$ and ${\rm ASTM}$ (for inside players) imply high values for ${\rm PMW}$.
        Thus, we will consider the following weights:
        \begin{itemize}
            \item Points guards: the weights for ${\rm EPts}$ y ${\rm ASTM}$ are set to $\omega_k=0.4$, while the rest of the criteria are set to $\omega_k=\dfrac{0.2}{5} = 0.04$.
            \item Centers: the weights for ${\rm DRM}$ y ${\rm ASTM}$ are set to $\omega_k=0.4$, while the rest of the criteria are set to $\omega_k=\dfrac{0.2}{5} = 0.04$.
        \end{itemize}
\end{itemize}
\item {\it Shape of the preference functions}: The Type V (V-Shape with indifference Criterion) preference function was consider to compare players for all the criteria. Such a function was chosen because the quantitative nature of the measure and also because its flexibility. We apply the method described in Section \ref{sec:promethee}, to determine the preference and indifference threshold, by considering $\alpha=25\%$ and $\beta=75\%$. This particular choice allows us, when comparing players, to get $25\%$ of the comparison resulting in high preference, and $25\%$ in indifferent. The remainder $50\%$ of the comparison result on preference degrees in $(0,1)$, which gives some flexibility to the decision process and allows to distinguish between players. The thresholds obtained with this choice to all the set of players and criteria are shown in Table \ref{table2}.

\end{itemize}

\begin{table}
    \centering
    \caption{Averages for each criteria when differentiating player's profiles. $p$-values of ANOVA test.}\label{anova}
    \begin{tabular}{cccccccc}
        \toprule
         & Total & PG & SG	& F & PF & C & ANOVA (p-value) \\
        \midrule
        PtsM &	0.371 &	0.346 &	0.372 &	0.353 &	0.390 &	0.396 &	0.081 \\
        DRM &	0.036 &	\textbf{0.008} &	-0.001 &	0.044 &	0.057 &	\textbf{0.077} &	$< 10^{-3}$ \\
        ORM &	0.145 &	0.166 &	0.094 &	0.100 &	0.147 &	0.222 &	$< 10^{-3}$ \\
        PMW &	0.216 &	0.126 &	0.368 &	0.401 &	-0.069 &	0.372 &	0.548 \\
        EPts &	45.993 &	42.279 &	41.462 &	43.508 &	47.592 &	56.510 &	$< 10^{-3}$ \\
        ASTM &	0.074 &	\textbf{0.116} &	0.077 &	0.067 &	0.055 &	\textbf{0.049} &	$< 10^{-3}$ \\
        PCS \% &	55.803 &	59.040 &	52.298 &	52.760 &	54.523 &	60.694 &	$< 10^{-3}$ \\
        \bottomrule
    \end{tabular}
\end{table}

\begin{table}
    \centering
    \caption{Minimum ($p$) and maximum ($q$) thresholds for each index depending on the player's position} \label{table2}
    {\small
    \begin{tabular}{ccccccc}
        \toprule
         &  &	Point guard &	Shooting guard &	Forward &	Power-forward &	Center \\
        \midrule
        \multirow{2}{*}{PtsM} & q &	0.045 &	0.043 &	0.041 &	0.035 &	0.033 \\
         & p &	0.164 &	0.158 &	0.163 &	0.131 &	0.143 \\\hline
        \multirow{2}{*}{DRM} &q &	0.019 &	0.019 &	0.029 &	0.026 &	0.029  \\
         & p &	0.067 &	0.071 &	0.104 &	0.097 &	0.104  \\\hline
        \multirow{2}{*}{ORM} & q &	0.033 &	0.035 &	0.028 &	0.035 &	0.027  \\
         & p &	0.136 &	0.143 &	0.111 &	0.152 &	0.096 \\\hline
        \multirow{2}{*}{EPts} & q &	3.042 &	2.875 &	2.354 &	2.675 &	1.855  \\
         & p &	12.091 &	10.935 &	9.088 &	10.241 &	7.195  \\\hline
        \multirow{2}{*}{ATSM} & q &	0.014 &	0.012 &	0.011 &	0.008 &	0.007  \\
         & p &	0.051 &	0.048 &	0.043 &	0.035 &	0.036  \\\hline
        \multirow{2}{*}{PCS \%} & q &	2.025 &	3.002 &	2.677 &	2.999 &	1.424  \\
         &p &	8.954 &	10.789 &	9.911 &	10.731 &	7.681  \\
         \bottomrule
    \end{tabular}
    }
\end{table}


\begin{table}
    \centering
    \caption{Correlation matrix of the criteria by player position (PG = Point Guards, SG = Shotting Guards, F = Forwards, PF= Power Forwards, C = Centers)} \label{table1}
    {\small
    \begin{tabular}{ccccccccc}
        \toprule
        & PtsM & DRM & ORM & PMW & EPts & ASTM & PCS \% \\\midrule
         & 1 & 0.13 & \textbf{0.502} & \textbf{0.503} & \textbf{0.604} & 0.084 & 0.24 & PG \\
         & 1 & -0.057 & 0.27 & 0.219 & \textbf{0.582} & -0.131 & 0.191 & SG \\
        PtsM & 1 & -0.032 & \textbf{0.501} & \textbf{0.348} & \textbf{0.592} & 0.085 & \textbf{0.438} & F \\
         & 1 & 0.279 & \textbf{0.417} & 0.115 & 0.276 & -0.182 & 0.237 & PF \\
         & 1 & 0.08 & \textbf{0.52} & 0.254 & -0.03 & 0.043 & -0.12 & C \\
         & & 1 & \textbf{0.52} & 0.122 & \textbf{0.47} & 0.179 & \textbf{0.514} & PG \\
         & & 1 & \textbf{0.346} & 0.007 & 0.291 & \textbf{0.461} & \textbf{0.501} & SG \\
        DRM  & & 1 & 0.337 & -0.11 & 0.032 & -0.266 & 0.131 & F \\
         & & 1 & 0.197 & \textbf{0.331} & 0.115 & \textbf{0.307} & \textbf{0.319} & PF \\
         & & 1 & \textbf{0.53} & \textbf{0.35} & \textbf{0.359} & 0.065 & \textbf{0.364} & C \\
         & & & 1 & \textbf{0.381} & \textbf{0.771} & \textbf{0.49} & \textbf{0.896} & PG \\
         & & & 1 & \textbf{0.432} & \textbf{0.727} & \textbf{0.412} & \textbf{0.875} & SG \\
        ORM  & & & 1 & \textbf{0.439} & \textbf{0.538} & 0.281 & \textbf{0.689} & F \\
         & & & 1 & 0.209 & \textbf{0.648} & 0.245 & \textbf{0.818} & PF \\
         & & & 1 & 0.245 & 0.32 & 0.171 & \textbf{0.603} & C \\
         & & & & 1 & \textbf{0.313} & \textbf{0.396} & \textbf{0.268} & PG \\
         & & & & 1 & \textbf{0.329} & 0.32 & \textbf{0.389} & SG \\
        PMW  & & & & 1 & \textbf{0.477} & \textbf{0.438} & \textbf{0.406} & F \\
         & & & & 1 & 0.257 & \textbf{0.386} & 0.286 & PF \\
         & & & & 1 & 0.18 & \textbf{0.532} & 0.241 & C \\
         & & & & & 1 & 0.142 & \textbf{0.743} & PG \\
         & & & & & 1 & 0.194 & \textbf{0.781} & SG \\
        EPts  & & & & & 1 & 0.082 & \textbf{0.734} & F \\
         & & & & & 1 & 0.071 & \textbf{0.746} & PF \\
         & & & & & 1 & \textbf{-0.354} & \textbf{0.503} & C \\
         & & & & &  & 1 & \textbf{0.564} & PG \\
         & & & & &  & 1 & \textbf{0.525} & SG \\
        ASTM  & & & & &  & 1 & 0.33 & F \\
         & & & & &  & 1 & \textbf{0.458} & PF \\
         & & & & &  & 1 & 0.186 & C \\
        \bottomrule
    \end{tabular}
    }
\end{table}

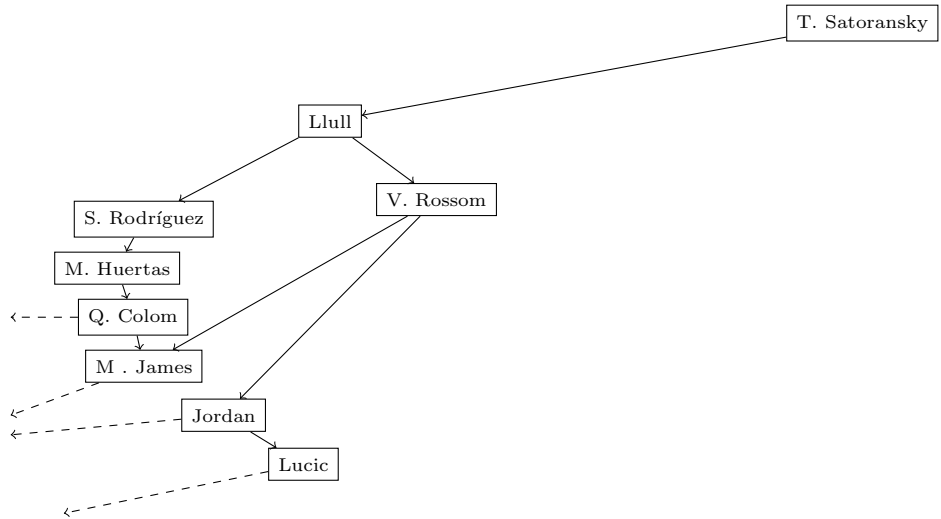
\begin{figure}
    \centering
    \begin{tikzpicture}[xscale=0.7, yscale=1.3]

\node[draw] (Satoransky) at (20,0) {\tiny T. Satoransky};
\node[draw] (Llull) at (10,-1) {\tiny Llull};
\node[draw] (Rossom) at (12,-1.8) {\tiny V. Rossom};
\node[draw] (Rodriguez) at (6.5,-2) {\tiny S. Rodr\'iguez};
\node[draw] (Huertas) at (6,-2.5) {\tiny M. Huertas};
\node[draw] (Colom) at (6.3,-3) {\tiny Q. Colom};
\node[draw] (James) at (6.5, -3.5) {\tiny M . James};
\node[draw] (Jordan) at (8,-4) {\tiny Jordan};
\node[draw] (Lucic) at (9.5,-4.5) {\tiny Lucic};

\draw[->] (Satoransky) --(Llull);
\draw[->] (Llull) --(Rossom);
\draw[->] (Llull) --(Rodriguez);
\draw[->] (Rossom) --(James);
\draw[->] (Rossom) --(Jordan);
\draw[->] (Rodriguez) --(Huertas);
\draw[->] (Huertas) -- (Colom);
\draw[->] (Colom) -- (James);
\draw[->, dashed] (James) -- (4, -4);
\draw[->, dashed] (Jordan) -- (4,-4.2);
\draw[->, dashed] (Colom) -- (4,-3);
\draw[->, dashed] (Lucic) -- (5,-5);
\draw[->] (Jordan) -- (Lucic);

\end{tikzpicture}
    \caption{Top part of the graph for points guard players under Scenario 1.\label{bases_sc1}}
\end{figure}

\begin{figure}
    \centering
    \begin{tikzpicture}[xscale=1, yscale=1.3]

\node[draw] (Satoransky) at (15,0) {\tiny T. Satoransky};
\node[draw] (Rossom) at (11,-1) {\tiny V. Rossom};
\node[draw] (Rodriguez) at (7.5,-2) {\tiny S. Rodr\'iguez};
\node[draw] (Huertas) at (5.5,-2) {\tiny M. Huertas};
\node[draw] (Llull) at (5,-2.5) {\tiny Llull};
\node[draw] (James) at (3.65, -3.5) {\tiny M . James};
\node[draw] (Jordan) at (5.15,-3.5) {\tiny Jordan};
\node[draw] (Markovic) at (9.5,-3) {\tiny Markovic};
\node[draw] (Schreiner) at (9,-3.5) {\tiny T. Schreiner};

\draw[->] (Satoransky) --(Rossom);
\draw[->] (Rossom) --(Rodriguez);
\draw[->] (Rossom) --(Huertas);
\draw[->] (Rodriguez) --(Llull);
\draw[->] (Huertas) --(Llull);
\draw[->] (Llull) --(James);
\draw[->] (Llull) --(Jordan);
\draw[->] (Llull) --(Markovic);
\draw[->, dashed] (Markovic) -- (13,-4);
\draw[->] (Markovic) -- (Schreiner);
\draw[->, dashed] (James) --(3,-4);
\draw[->, dashed] (Jordan) --(3.2,-4);
\draw[->, dashed] (Schreiner) --(8,-4);

\end{tikzpicture}
    \caption{Top part of the graph for points guard players under Scenario 2.\label{bases_sc2}}
\end{figure}
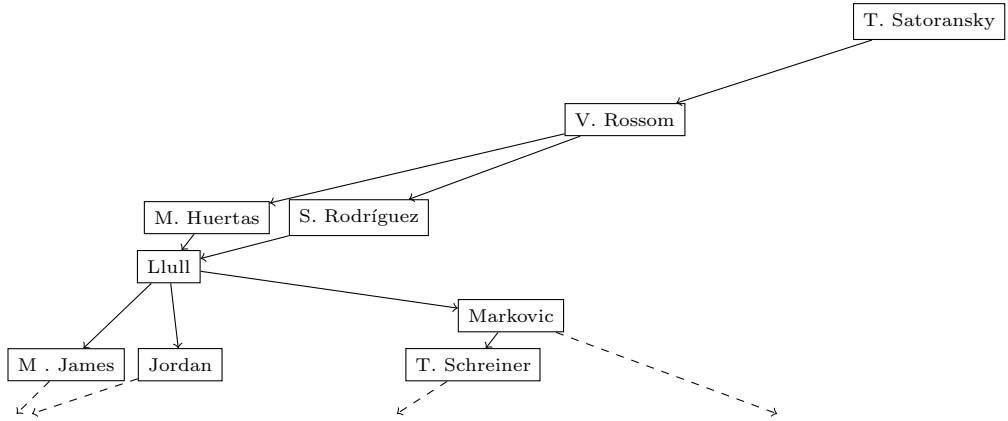

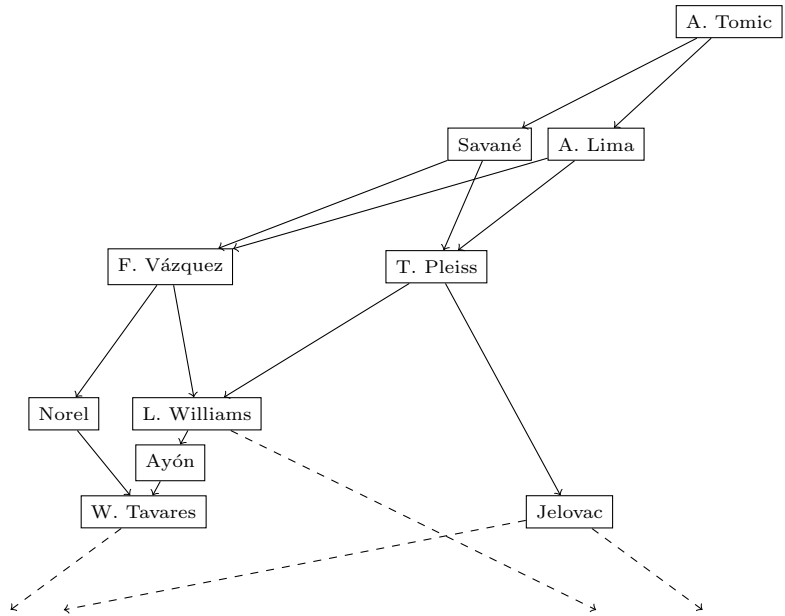
\begin{figure}
    \centering
    \begin{tikzpicture}[xscale=0.7, yscale=1.3]

\node[draw] (Tomic) at (17.5,0) {\tiny A. Tomic};
\node[draw] (Savane) at (13,-1.25) {\tiny Savan\'e};
\node[draw] (Lima) at (15,-1.25) {\tiny A. Lima};
\node[draw] (Pleiss) at (12,-2.5) {\tiny T. Pleiss};
\node[draw] (Vazquez) at (7,-2.5) {\tiny F. V\'azquez};
\node[draw] (Norel) at (5,-4) {\tiny Norel};
\node[draw] (Williams) at (7.5, -4) {\tiny L. Williams};
\node[draw] (Ayon) at (7,-4.5) {\tiny Ay\'on};
\node[draw] (Tavares) at (6.5,-5) {\tiny W. Tavares};
\node[draw] (Jelovac) at (14.5,-5) {\tiny Jelovac};

\draw[->] (Tomic) --(Savane);
\draw[->] (Tomic) --(Lima);
\draw[->] (Savane) --(Pleiss);
\draw[->] (Lima) --(Pleiss);
\draw[->] (Savane) --(Vazquez);
\draw[->] (Lima) --(Vazquez);
\draw[->] (Pleiss) --(Williams);
\draw[->] (Pleiss) --(Jelovac);
\draw[->] (Norel) --(Tavares);
\draw[->] (Vazquez) --(Norel);
\draw[->] (Vazquez) --(Williams);
\draw[->] (Williams) --(Ayon);
\draw[->] (Ayon) --(Tavares);
\draw[->, dashed] (Tavares)-- (4,-6);
\draw[->, dashed] (Jelovac)-- (5,-6);
\draw[->, dashed] (Jelovac)-- (17,-6);
\draw[->, dashed] (Williams)-- (15,-6);

\end{tikzpicture}
    \caption{Top part of the graph for center players under Scenario 1.\label{pivots_sc1}}
\end{figure}

\begin{figure}
    \centering
    \begin{tikzpicture}[xscale=0.7, yscale=1.4]

\node[draw] (Tomic) at (17.5,0) {\tiny A. Tomic};
\node[draw] (Ayon) at (14.5,-1) {\tiny Ay\'on};
\node[draw] (Lima) at (15,-1.5) {\tiny A. Lima};
\node[draw] (Savane) at (11,-2) {\tiny Savan\'e};
\node[draw] (Tavares) at (9,-3) {\tiny W. Tavares};
\node[draw] (Llompart) at (15,-3) {\tiny Llompart};
\node[draw] (Loncar) at (18.5,-5) {\tiny Loncar};
\node[draw] (Bourouisis) at (8.5,-3.5) {\tiny Bourouisis};

\node[draw] (Begic) at (8.25,-4.5) {\tiny M. Begic};
\node[draw] (Tillie) at (10, -4.75) {\tiny K. Tillie};

\draw[->] (Tomic) --(Ayon);
\draw[->] (Ayon) --(Lima);
\draw[->] (Lima) --(Savane);
\draw[->] (Savane) --(Llompart);
\draw[->] (Savane) --(Tavares);
\draw[->] (Tavares) --(Bourouisis);
\draw[->] (Bourouisis)--(Begic);
\draw[->] (Bourouisis)--(Tillie);
\draw[->] (Llompart)--(Tillie);
\draw[->] (Llompart)--(Loncar);

\draw[->, dashed] (Tavares)-- (4,-4);
\draw[->, dashed] (Bourouisis)-- (5,-4.36);
\draw[->, dashed] (Begic)-- (4.5,-5.02);
\draw[->, dashed] (Tillie)-- (4,-5.3);
\draw[->, dashed] (Loncar)-- (4,-5.7);

\end{tikzpicture}
    \caption{Top part of the graph for center players under Scenario 2.\label{pivots_sc2}}
\end{figure}
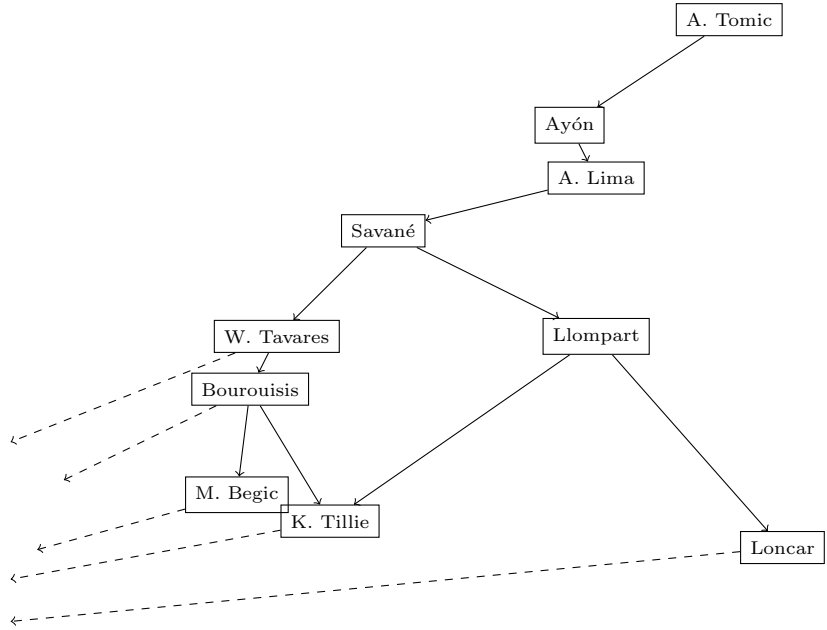

The flows obtained are for each of the two scenarios and each of the two profiles of players are shown in Table \ref{table:fbases} (Point Guards) and Table \ref{table:fpivots} (Centers). The PROMETHEE I graphs which can be built using the preference flows are drawn in Figures \ref{bases_sc1}-\ref{pivots_sc2}.

As can be observed from Figures \ref{bases_sc1} and \ref{bases_sc2}, the most preferable point guard player under the two scenarios is Satoranski. This player obtained the maximum values for three of the six considered criteria (ORM, EPts and PSC\%) and his results are among the best 25\% for the the remainder criteria, justifying the top 1 ranking in this profile. On the other hand, Llull obtained a clear 2nd position under Scenario 1, while this position was obtained by V. Rossom under Scenario 2. This \textit{change} of positions can be explained after a simple analysis of the absolute data obtained by each player:

\begin{center}
    \begin{tabular}{rcccccc}
        \toprule
        & PtsE & DRM & ORM & EPts & ASTM & PCS\%\\
        \midrule
        Llull & 0.49 & 0.06 & 0.29 & {\bf 49.08} & {\bf 0.14} &  65.1\\
        V. Rossom & 0.45 & -0.02 & 0.26 & {\bf 51.2} & {\bf 0.17} & { 57.07}\\
        \bottomrule
    \end{tabular}
\end{center}

\noindent where it can be observed that the criteria whose weights have been increased in Scenario 2 (EPts and ASTM) where those in which Rossom got better results that Llull. A similar consideration can be established for Markovic and Schreiner whose results in ASTM where the highest, increasing its preference under Scenario 2 with respect to those obtained under Scenario 1. Observe that this analysis which seems to be straightforward when comparing two players, can be cumbersome when analyzing the whole sample (42 players).

With respect to the Center players (Figures \ref{pivots_sc1} and \ref{pivots_sc2}), Tomic was ranked in the first position under the two scenarios. Some remarkable considerations when comparing the rankings obtained with the two scenarios are the following:
\begin{itemize}
\item Ay\'on considerably increases its ranking position under Scenario 2 w.r.t the one obtained under Scenario 1. It follows since the criteria whose  weights were increased (DRM and ASTM) under Scenario 2, were those for which Ay\'on was among the top 5 players with best results.
\item T. Pleiss was ranked among the best players under Scenario 1, since it is the one with best result in EPts, decreases significatively its position under Scenario 2, mainly because the weight for EPts was 0.04 under this scenario.
\end{itemize}

Observe that in some cases the methodology PROMETHEE I does not give us a total order over the players. To obtain an strict ranking, we use the methodology PROMETHEE II, which aggregate the positive and negative preference flows into a single measure (net preference flow). Hence, although in some case only little differences, this approach allows us to strictly rank the players.

For instance, under Scenario 1, V. Rossom and S. Rodr\'iguez (Point guard players) or Lima and Savan\'e (Center players) are among those cases. PROMETHEE I is not able to determine which dominates (the positive flow is greater and the negative flow is smaller for a player) but PROMETHEE II states that V. Rossom and Lima obtained (slightly) larger net flow values than S. Rodr\'iguez and Savan\'e, respectively (see preference flows for center and guard players in Tables \ref{table:fbases} and \ref{table:fpivots}).

%
%

\begin{table}
  \centering
  \caption{PROMETHEE flows for point guard players.\label{table:fbases}}%
  {\small
  \begin{tabular}{lcccclcccc}
  \toprule
    \multicolumn{4}{c}{Scenario 1} & & \multicolumn{4}{c}{Scenario 2}\\
    \textbf{Players} & $\Phi$ & $\Phi^+$ &$\Phi^-$ &       & \textbf{Players} & $\Phi$ & $\Phi^+$ &$\Phi^-$ \\
    \midrule
    \textbf{T. Satoransky} & 0.7222 & 0.7343 & 0.0120 &       & \textbf{T. Satoransky} & 0.7246 & 0.7378 & 0.0132 \\
    \textbf{Llull} & 0.4918 & 0.5178 & 0.0260 &       & \textbf{V. Rossom} & 0.5644 & 0.5950 & 0.0306 \\
    \textbf{V. Rossom} & 0.4176 & 0.4932 & 0.0756 &       & \textbf{S. Rodr\'iguez} & 0.4288 & 0.4647 & 0.0359 \\
    \textbf{S. Rodr\'iguez} & 0.4069 & 0.4425 & 0.0356 &       & \textbf{M. Huertas} & 0.4281 & 0.4629 & 0.0348 \\
    \textbf{M. Huertas} & 0.3622 & 0.4051 & 0.0429 &       & \textbf{Llull} & 0.3836 & 0.4319 & 0.0483 \\
    \textbf{Q. Colom} & 0.3108 & 0.3822 & 0.0714 &       & \textbf{Markovic} & 0.2915 & 0.4272 & 0.1357 \\
    \textbf{M. James} & 0.2655 & 0.3631 & 0.0976 &       & \textbf{Jordan} & 0.2742 & 0.3636 & 0.0893 \\
    \textbf{San Miguel} & 0.2476 & 0.3242 & 0.0767 &       & \textbf{M. James} & 0.2742 & 0.3518 & 0.0777 \\
    \textbf{Jordan} & 0.2407 & 0.3717 & 0.1310 &       & \textbf{T. Schreiner} & 0.1981 & 0.3736 & 0.1756 \\
    \textbf{Granger} & 0.2358 & 0.3245 & 0.0887 &       & \textbf{T. Bellas} & 0.1920 & 0.2998 & 0.1078 \\
    \textbf{T. Bellas} & 0.2336 & 0.3115 & 0.0779 &       & \textbf{V. Sada} & 0.1886 & 0.3638 & 0.1753 \\
    \textbf{A. Oliver} & 0.2275 & 0.3172 & 0.0896 &       & \textbf{San Miguel} & 0.1861 & 0.2859 & 0.0998 \\
    \textbf{Lucic} & 0.1936 & 0.3615 & 0.1680 &       & \textbf{A. Oliver} & 0.1627 & 0.2663 & 0.1037 \\
    \textbf{R. Neto} & 0.1522 & 0.2652 & 0.1130 &       & \textbf{Q. Colom} & 0.1613 & 0.2696 & 0.1084 \\
    \textbf{R. Luz} & 0.1056 & 0.2573 & 0.1517 &       & \textbf{R. Neto} & 0.1335 & 0.2680 & 0.1345 \\
    \textbf{C. Cabezas} & 0.0912 & 0.2342 & 0.1429 &       & \textbf{D. Perkins} & 0.1108 & 0.3212 & 0.2104 \\
    \textbf{V. Sada} & 0.0819 & 0.2795 & 0.1976 &       & \textbf{D. White} & 0.0970 & 0.2590 & 0.1620 \\
    \textbf{Fern\'andez} & 0.0788 & 0.2201 & 0.1413 &       & \textbf{R. Luz} & 0.0908 & 0.2533 & 0.1625 \\
    \textbf{Markovic} & 0.0737 & 0.2708 & 0.1971 &       & \textbf{C. Bivi\`a} & 0.0901 & 0.2240 & 0.1339 \\
    \textbf{Woodside} & 0.0475 & 0.1998 & 0.1523 &       & \textbf{Granger} & 0.0875 & 0.2317 & 0.1442 \\
    \textbf{T. Schreiner} & 0.0451 & 0.2366 & 0.1915 &       & \textbf{Lucic} & 0.0748 & 0.3959 & 0.3211 \\
    \textbf{C. Bivi\`a} & 0.0395 & 0.1915 & 0.1520 &       & \textbf{Woodside} & 0.0666 & 0.2112 & 0.1446 \\
    \textbf{D. White} & 0.0308 & 0.2328 & 0.2020 &       & \textbf{Vives} & 0.0484 & 0.2182 & 0.1698 \\
    \textbf{Hannah} & 0.0268 & 0.2584 & 0.2316 &       & \textbf{R. L\'opez} & -0.0058 & 0.1849 & 0.1906 \\
    \textbf{Vives} & 0.0128 & 0.1886 & 0.1758 &       & \textbf{Fern\'andez} & -0.0080 & 0.1814 & 0.1894 \\
    \textbf{R. L\'opez} & -0.0456 & 0.1771 & 0.2227 &       & \textbf{C. Cabezas} & -0.0087 & 0.2057 & 0.2144 \\
    \textbf{T. Heurtel} & -0.0929 & 0.1593 & 0.2522 &       & \textbf{Hannah} & -0.0994 & 0.1999 & 0.2993 \\
    \textbf{D. Perkins} & -0.1087 & 0.2058 & 0.3145 &       & \textbf{Mallet} & -0.1052 & 0.1696 & 0.2748 \\
    \textbf{Radicevic} & -0.1241 & 0.1775 & 0.3016 &       & \textbf{Radicevic} & -0.1363 & 0.1486 & 0.2849 \\
    \textbf{A. Hern\'andez} & -0.1259 & 0.1710 & 0.2970 &       & \textbf{P. Pozas} & -0.1626 & 0.2455 & 0.4081 \\
    \textbf{Mallet} & -0.1282 & 0.1857 & 0.3139 &       & \textbf{T. Heurtel} & -0.2013 & 0.1160 & 0.3173 \\
    \textbf{Lisch} & -0.1844 & 0.1030 & 0.2874 &       & \textbf{A. D\'iaz} & -0.2045 & 0.1624 & 0.3669 \\
    \textbf{A. D\'iaz} & -0.2225 & 0.1567 & 0.3792 &       & \textbf{Franch} & -0.2114 & 0.1068 & 0.3182 \\
    \textbf{P. Pozas} & -0.2621 & 0.1413 & 0.4034 &       & \textbf{Lisch} & -0.2608 & 0.0878 & 0.3486 \\
    \textbf{Franch} & -0.2838 & 0.0887 & 0.3725 &       & \textbf{A. Hern\'andez} & -0.2934 & 0.1220 & 0.4154 \\
    \textbf{D. Adams} & -0.3178 & 0.1510 & 0.4688 &       & \textbf{Salgado} & -0.3378 & 0.0947 & 0.4326 \\
    \textbf{D. P\'erez} & -0.3194 & 0.0841 & 0.4035 &       & \textbf{D. P\'erez} & -0.3531 & 0.0771 & 0.4302 \\
    \textbf{Salgado} & -0.3474 & 0.0913 & 0.4387 &       & \textbf{D. Adams} & -0.3874 & 0.0891 & 0.4765 \\
    \textbf{R. Grimau} & -0.3683 & 0.0814 & 0.4497 &       & \textbf{R. Grimau} & -0.4150 & 0.0785 & 0.4935 \\
    \textbf{J. Mayo} & -0.6521 & 0.0364 & 0.6885 &       & \textbf{A. Rodr\'iguez} & -0.5272 & 0.0889 & 0.6161 \\
    \textbf{A. Rodr\'iguez} & -0.7129 & 0.0411 & 0.7541 &       & \textbf{J. Mayo} & -0.6491 & 0.0320 & 0.6811 \\
    \textbf{R. Huertas} & -0.8457 & 0.0020 & 0.8477 &       & \textbf{R. Huertas} & -0.8907 & 0.0044 & 0.8951 \\
    \bottomrule
    \end{tabular}%
    }
\end{table}

\begin{table}
  \centering
  \caption{PROMETHEE flows for center players.\label{table:fpivots}}%
  {\small
  \begin{tabular}{lcccclcccc}
    \toprule
    \multicolumn{4}{c}{Scenario 1} & & \multicolumn{4}{c}{Scenario 2}\\
    \textbf{Players} & $\Phi$ & $\Phi^+$ &$\Phi^-$ &       & \textbf{Players} & $\Phi$ & $\Phi^+$ &$\Phi^-$ \\
    \midrule
    \textbf{A. Tomic} & 0.5818 & 0.6095 & 0.0277 &       & \textbf{A. Tomic} & 0.6721 & 0.6855 & 0.0134 \\
    \textbf{A. Lima} & 0.4970 & 0.5413 & 0.0443 &       & \textbf{Ay\'on} & 0.5696 & 0.6059 & 0.0363 \\
    \textbf{Savan\'e} & 0.4963 & 0.5345 & 0.0383 &       & \textbf{A. Lima} & 0.5488 & 0.6031 & 0.0543 \\
    \textbf{T. Pleiss} & 0.3120 & 0.4238 & 0.1117 &       & \textbf{Savan\'e} & 0.4504 & 0.5120 & 0.0616 \\
    \textbf{F. V\'azquez} & 0.3039 & 0.3843 & 0.0804 &       & \textbf{W. Tavares} & 0.2707 & 0.3932 & 0.1225 \\
    \textbf{L. Williams} & 0.2466 & 0.3658 & 0.1193 &       & \textbf{Llompart} & 0.2653 & 0.4653 & 0.2000 \\
    \textbf{Norel} & 0.2375 & 0.3327 & 0.0952 &       & \textbf{Bourouisis} & 0.2099 & 0.3571 & 0.1473 \\
    \textbf{Ay\'on} & 0.2221 & 0.3454 & 0.1233 &       & \textbf{M. Begic} & 0.1406 & 0.3180 & 0.1774 \\
    \textbf{W. Tavares} & 0.1986 & 0.3307 & 0.1321 &       & \textbf{C. Iverson} & 0.1333 & 0.2577 & 0.1244 \\
    \textbf{Jelovac} & 0.1699 & 0.3746 & 0.2047 &       & \textbf{K. Tillie} & 0.1248 & 0.3313 & 0.2065 \\
    \textbf{M. Begic} & 0.1128 & 0.2689 & 0.1561 &       & \textbf{W. Hern\'angomez} & 0.1205 & 0.2584 & 0.1379 \\
    \textbf{C. Iverson} & 0.0967 & 0.2464 & 0.1497 &       & \textbf{F. V\'azquez} & 0.0828 & 0.2192 & 0.1364 \\
    \textbf{Loncar} & 0.0821 & 0.3538 & 0.2717 &       & \textbf{Jelovac} & 0.0801 & 0.2484 & 0.1683 \\
    \textbf{J. Akindele} & 0.0423 & 0.2072 & 0.1649 &       & \textbf{Norel} & 0.0781 & 0.2194 & 0.1413 \\
    \textbf{W. Hern\'angomez} & 0.0382 & 0.2267 & 0.1885 &       & \textbf{Loncar} & 0.0416 & 0.3995 & 0.3580 \\
    \textbf{N. Jawai} & -0.0030 & 0.2267 & 0.2297 &       & \textbf{Sekulic} & 0.0316 & 0.2272 & 0.1955 \\
    \textbf{Sekulic} & -0.0368 & 0.2020 & 0.2388 &       & \textbf{J. Akindele} & 0.0049 & 0.1906 & 0.1857 \\
    \textbf{H. Rizvic} & -0.0986 & 0.1847 & 0.2834 &       & \textbf{L. Williams} & -0.0481 & 0.1792 & 0.2274 \\
    \textbf{Llompart} & -0.1121 & 0.2822 & 0.3944 &       & \textbf{T. Pleiss} & -0.0634 & 0.1581 & 0.2215 \\
    \textbf{Golubovic} & -0.1144 & 0.1306 & 0.2450 &       & \textbf{Golubovic} & -0.1037 & 0.1349 & 0.2385 \\
    \textbf{J. Triguero} & -0.1352 & 0.1751 & 0.3103 &       & \textbf{H. Rizvic} & -0.1048 & 0.1224 & 0.2272 \\
    \textbf{Lampropoulos} & -0.1755 & 0.0982 & 0.2736 &       & \textbf{Balvin} & -0.1169 & 0.1576 & 0.2745 \\
    \textbf{Bourouisis} & -0.1787 & 0.1649 & 0.3435 &       & \textbf{Lampropoulos} & -0.2137 & 0.0740 & 0.2876 \\
    \textbf{K. Tillie} & -0.2395 & 0.1567 & 0.3962 &       & \textbf{M. Diagne} & -0.2539 & 0.0538 & 0.3077 \\
    \textbf{Balvin} & -0.2512 & 0.1251 & 0.3763 &       & \textbf{J. Triguero} & -0.2700 & 0.0754 & 0.3455 \\
    \textbf{Doblas} & -0.2825 & 0.0579 & 0.3405 &       & \textbf{G. Bogris} & -0.3043 & 0.1040 & 0.4083 \\
    \textbf{G. Bogris} & -0.2961 & 0.0930 & 0.3892 &       & \textbf{N. Jawai} & -0.3060 & 0.0656 & 0.3716 \\
    \textbf{M. Diagne} & -0.3065 & 0.0610 & 0.3675 &       & \textbf{Doblas} & -0.3494 & 0.0252 & 0.3746 \\
    \textbf{Katic} & -0.3324 & 0.0828 & 0.4151 &       & \textbf{Miralles} & -0.3739 & 0.0221 & 0.3960 \\
    \textbf{Miralles} & -0.3386 & 0.0556 & 0.3942 &       & \textbf{D. Miller} & -0.3868 & 0.0369 & 0.4237 \\
    \textbf{Lishchuk} & -0.3397 & 0.0620 & 0.4017 &       & \textbf{Katic} & -0.4429 & 0.0342 & 0.4771 \\
    \textbf{D. Miller} & -0.3970 & 0.0718 & 0.4688 &       & \textbf{Lishchuk} & -0.4873 & 0.0132 & 0.5005 \\
    \bottomrule
    \end{tabular}%
    }
\end{table}

\section{Conclusions}
    \label{sec:conclusions}

In this paper we present a flexible quantitative tool that allows us to aggregate some of the well-known indices to measure the basketball players' skills, to provide either a partial (but richer) or total ranking on the players. The methodology is based on the use of the multicriteria method PROMETHEE. The PROMETHEE methodologies have the advantage that they do not consider the single player information, but the role of the player in the whole context, using the differences of the results between each two players instead of the single absolute value of the criterion for each player. Although this method is based on the interaction between the provided information and the experience provided by the decision maker, we propose an alternative approach to determine the values of some of the parameters which are needed for the use of the PROMETHEE methods. 

One of the highlights of the method is that the results can be presented using a graph, and that can be easily interpreted by practitioners non familiar with the mathematical tool. That is, this work shows that performance data can be used by agents, managers or coaches, by adapting their needs and choices on what type of players they have or need in different seasons. That is, by using our approach, coaches are able to obtain very useful information to help them in the decision of signing a player up, or in case of players and agents, such an information allows them to contextualize the player's performance, helping them to improve their efficiency or to get better contracts. 

It is left for further research the incorporation of economic (budget constraints) or qualitative (like leadership or discipline) information of the players, either as one of the criterion or using a post-process optimization tools that, based on the PROMETHEE results allows us to adapt the ``best'' players to the team budget. Also, in \cite{corrente}, the combination of Stochastic Multicriteria Acceptability Analysis (SMAA) and the PROMETHEE method to determine the adequate weights for the decision. The application of such an approach will be the topic of a forthcoming paper.

\section*{Acknowledgments}

The first and second authors were partially supported by project PP2016-PIP06 (Universidad de Granada) and the research group SEJ-534 (Junta de Andaluc\'ia). The first author was also partially supported by project MTM2016-74983-C2-1-R (MINECO, Spain).

\end{document}